\documentclass[10pt]{amsart}%
\usepackage[utf8]{inputenc}%
\usepackage[english]{babel}%
\usepackage{graphicx}%
\usepackage{amsmath}%
\usepackage{amssymb}%
\usepackage{amsthm}%
\usepackage{amsfonts}%
\usepackage{hyperref}
\newcommand{\Kb}{\mathbb{K}}
\newcommand{\m}{\mathbf{m}}
\newcommand{\al}{\alpha}
\newcommand{\be}{\beta}
\newcommand{\te}{\theta}
\newcommand{\la}{\lambda}
\newcommand{\La}{\Lambda}
\newcommand{\N}{\mathbb{Z}_{>0}}
\newcommand{\Nn}{\mathbb{Z}_{\ge0}}
\renewcommand{\ge}{\geqslant}
\renewcommand{\le}{\leqslant}
\newcommand{\da}{\downarrow}
\newcommand{\ua}{\uparrow}
\newcommand{\krbig}[2]{\big[{#1}:{#2}\big]}
\newcommand{\kr}[2]{[{#1}:{#2}]}
\newcommand{\kra}[2]{\kappa({#2},{#1})}
\newcommand{\Lao}{\La^\circ}
\newcommand{\sk}{\overline{\nabla}_\infty}
\newcommand{\Fc}{\mathcal{F}}
\newcommand{\pd}[2]{\mathrm{PD}(#1,#2)}
\newcommand{\x}[2]{\mathbf{X}_{#1,#2}}
\newcommand{\con}{\mathrm{r}}
\newtheorem{thm}{Theorem}[section]
\newtheorem{prop}[thm]{Proposition}
\newtheorem{lemma}[thm]{Lemma}
\newtheorem{df}[thm]{Definition}
\newtheorem{note}[thm]{Remark}
\begin{document}

\title[Diffusion processes in the 
Kingman simplex]{A two-parameter family of infinite-dimensional diffusions in the Kingman simplex}
\author{Leonid Petrov}
\maketitle

\section*{Introduction}
The main result of the present paper is to construct a two-parameter
family of Markov processes 
$\x{\al}{\te}(t)$ in the infinite-dimensional Kingman
simplex
\begin{equation*}
	\sk=\biggl\{x=(x_1,x_2,\dots)\in[0,1]^\infty\colon x_1\ge x_2\ge\dots\ge0,\
	{\displaystyle \sum\limits_{i=1}^\infty x_i\le 1}\biggr\}.
\end{equation*}
In the topology of coordinatewise convergence $\sk$ is a compact, metrizable
and separable space. Denote by $C(\sk)$ the algebra of real continuous
functions on $\sk$ with pointwise operations and the supremum norm.

In $C(\sk)$ there is a distinguished dense subspace $\Fc:=\mathbb{R}\left[ q_1,q_2,\dots \right]$
generated (as a commutative unital algebra) by algebraically independent
continuous functions 
$q_k(x):=\sum_{i=1}^{\infty}x_i^{k+1}$, 
$k=1,2,\dots$, $x\in\sk$.

For each $0\le\al<1$ and $\te>-\al$ we define an operator $A\colon\Fc\to\Fc$
which can be written as a formal differential operator of second order
with respect to the generators of the algebra $\Fc$:\footnote{Here, by agreement, $q_0=1$.}
\begin{equation}\label{f0.1}
	\begin{array}{lcl}
	A&=&\displaystyle 
	\sum_{i,j=1}^{\infty}(i+1)(j+1)(q_{i+j}-q_iq_j)\frac{\partial^2}{\partial q_i\partial q_j}
	\\&&\displaystyle\qquad+
	\sum_{i=1}^{\infty}\left[ -(i+1)(i+\te)q_i+(i+1)(i-\al)q_{i-1} \right]\frac{\partial}{\partial q_i},
	\end{array}
\end{equation}
or (subject to certain restrictions, see Remarks \ref{p5.2} and \ref{p5.3} below) as a
differential operator in natural coordinates:
\begin{equation}\label{f0.2}
	A=\sum_{i=1}^{\infty}x_i\frac{\partial^2}{\partial x_i^2}-
	\sum_{i,j=1}^{\infty}x_ix_j\frac{\partial^2}{\partial x_i\partial x_j}-
	\sum_{i=1}^{\infty}(\te x_i+\al)\frac{\partial}{\partial x_i}.
\end{equation}
It is worth noting that the operator $A$ is defined only on $\Fc$, and 
a direct application
of the right-hand sides of 
(\ref{f0.1}) and (\ref{f0.2}) to other functions
requires caution.

\smallskip
{\bf{}Main theorem\/}.
{\em{}	{\rm{}(1)\/} 
	The operator $A$ is closable in $C(\sk)$;

	{\rm{}(2)\/} 
	Its closure $\overline A$ generates
	a diffusion process 
	$\{\x{\al}{\te}(t)\}_{t\ge0}$ in $\sk$,
	that is, a strong Markov process with continuous sample paths;

	{\rm{}(3)\/} 
	The two-parameter Poisson-Dirichlet measure 
	$\pd{\al}{\te}$
	on $\sk$ 
	is a unique invariant probability measure for
	$\x{\al}{\te}(t)$. 
	The process is reversible with respect to that measure;

	{\rm{}(4)\/} 
	The spectrum of the generator $\overline A$ 
	is described explicitly in \S\ref{s4.2} below.
	Due to existence of a spectral gap,
	the process 
	$\x{\al}{\te}(t)$
	is ergodic with respect to the measure 
	$\pd{\al}{\te}$\/}.

\smallskip
Fix arbitrary $0\le\al<1$ and $\te>-\al$. 
We construct the process in the simplex $\sk$
as a limit for a sequence of finite Markov chains 
with growing number of states.
The state space for $n$th chain
is the set of all partitions of $n$.
Denote this set by $\Kb_n$.
Each $\Kb_n$ is equipped with a certain probability measure $M_n$
(also depending on $\al$ and $\te$)
which is explicitly written out.
The system $\{M_n\}$ is a partition structure.\footnote{This
concept was introduced in \cite{Ki3}.
In a more general setting also
the name ``coherent system of random partitions''
is used (see, e.g., \cite{KOO,BO}).}
A bijection between partition structures and probability
measures on $\sk$ was established in \cite{Ki3}.
The Poisson-Dirichlet measure $\pd{\al}{\te}$ on $\sk$
corresponds to the two-parameter 
Ewens-Pitman's partition struc\-ture $\{M_n\}$.
In the case $\al=0$ this
partition structure
was introduced in \cite{Ew}, 
the two-parameter generalization is
due to Pitman \cite{Pi1}.
A special case $\al=0$, $\te=1$ of
$\{M_n\}$ was considered in \cite{VS}
in connection with 
limit behaviour of certain functionals
on the symmetric group $\mathfrak{S}_n$ as $n\to\infty$.
Also in \cite{VS} the measure $\pd{0}{1}$ was studied.
Two-parameter Ewens-Pitman's partition structures
were studied in
\cite{Pi2,Pi-Rand,Pi3,GP,PY,DGP}
and many other works. 
See also \cite{Pi3} for more bibliography.

It should be noted that the first example of a Markov process in
$\sk$ having the two-parameter Poisson-Dirichlet measure as a unique
invariant symmetrizing probability measure was constructed in 
\cite{Be}.
Diffusions
in the infinite-dimensional unordered simplex preserving the 
GEM-distribution (Griffths-Engen-McCloskey, see \cite{Pi3}) 
are considered in
a recent paper \cite{FW}. 
The GEM-distribution maps into
the Poisson-Dirichlet 
measure if we reorder the 
coordinates in a descending order.

The case $\al=0$ 
is relatively well understood
and
has applications 
to population genetics.
The measure $\pd{0}{\te}$ appeared in \cite{Ki1}.
The
process $\x{0}{\te}(t)$ was const\-ructed in 
\cite{EK1} as a limit for diffusion processes
in finite-dimensional 
simplexes as their dimension grows. 
In \S\ref{s5.4} we show that 
these finite-dimensional diffusions
arise as a particular
case of the process $\x{\al}{\te}(t)$
when $\al=-\be<0$ and $\te=N\be$, $N=2,3,\dots$.
Note also that when $\al=0$ the operator $A$ is equal to the operator
\cite[(2.10)]{EK1} multiplied by two.

In \cite{EK1} it was also proved that 
$\x{0}{\te}(t)$
is a limit process for a certain
Moran-type population model (about this model see also \cite{Wa1} and
\cite[Model II]{KMG}). 
This Moran-type model is a sequence of finite
Markov chains on partitions. They differ from the ones considered in
the present paper but have the same limit.

The process $\x{0}{\te}(t)$ (when $\al=0$)
is called the {\em{}infinitely
many neutral alleles diffusion model\/}.
It was also studied in 
\cite{Et,Schm} and other works.
The results of the present paper
extend some of the results
of \cite{EK1} and \cite{Et}
to the case $\al\ne0$.
However, it seems that in the general
case there is no connection between the model $\x{\al}{\te}(t)$ and population
genetics. Namely, both finite Markov chains and the process in $\sk$
have no known interpretation in terms of population genetics.

In the present paper the construction of 
the process $\x{\al}{\te}(t)$ in $\sk$ involves
only partition structures and up/down Markov chains. 
The up/down chains first appeared in \cite{Fu1}.
They were also studied in the papers \cite{BO,Fu2}.
The setting of the problem studied in the present paper 
and the general approach to it were
inspired by \cite{BO}. 
However, the concrete computations here
are performed in a different way. Monomial symmetric functions in
the present paper play the same role as Schur functions in \cite{BO}. 
The former are simpler than the latter, and the final result
is achieved by simpler means.
In addition to the results similar to \cite{BO}
one can
compute the operator $A$ in natural coordinates (\ref{f0.2}) and describe the
process $\x{\al}{\te}(t)$ when $\al=-\be<0$ and $\te=N\be$. 
Let us describe the organization of the paper.

In \S\ref{s1.1} 
we recall some notation concerning partitions.
In \S\ref{s1.2} the definitions
of partition structures and up/down Markov chains related to them are recalled.
In \S\ref{s1.3} we give the definition 
of Ewens-Pitman's partition structures. 
We also recall some of their properties.

In \S\ref{s2} 
we deal with some properties of 
symmetric functions
in the coordinates $\left( x_1,x_2,\dots \right)$
of a point $x\in\sk$.
These symmetric functions form the algebra $\Fc$ defined above.
In these terms we formulate the Kingman
theorem about the one-to-one
correspondence between partition structures and 
probability measures on $\sk$.

In \S\ref{s3} 
we obtain an explicit 
expression for the action
of operators $T_n$
(each $T_n$ is a transition operator of the 
$n$th up/down Markov chain)
on symmetric functions in the coordinates
$\la_1,\la_2,\dots$ of a partition $\la\in\Kb_n$.

In \S\ref{s4.1} 
we perform a limit transition from the up/down Markov chains
to the process $\x{\al}{\te}(t)$ on $\sk$.
First, we use the connection of symmetric functions in the 
coordinates of a partition with symmetric functions on $\sk$
to explicitly compute the limit of the operators $n^2(T_n-{\bf1})$.
This limit is an operator $A\colon\Fc\to\Fc$ described above.
After that we use some general results from the book \cite{EK2}
to establish the convergence of discrete semigroups 
$\left\{ {\bf1},T_n,T_n^2,\dots \right\}$
to the continuous semigroup $\{T(t)\}_{t\ge0}$
of 
the process
$\x{\al}{\te}(t)$.
This semigroup is generated by the closure 
$\overline A$
of the operator $A\colon\Fc\to\Fc$.
In \S\ref{s4.2} we formulate the remaining results
of the Main Theorem except the result about the continuity
of the process $\x{\al}{\te}(t)$'s sample paths.

In \S\ref{s5.1} 
we derive (\ref{f0.2}).
We also give explanations why the RHS of this formula should not be 
understood literally.
In \S\ref{s5.2} 
we derive the formula (\ref{f0.1}) for the operator $A$.
Using this formula it is possible
to prove the continuity of the sample paths of the process in $\sk$.
In \S\ref{s5.4} 
we deal with the case of 
degenerate values of the parameters, namely,
$\al=-\be$, $\te=N\be$.
In this case as a limit of up/down Markov chains
instead of the diffusion in $\sk$
we obtain a diffusion in the $(N-1)$-dimensional
simplex. This diffusion coincides with the one
studied in \cite{EK1}.

The author is very grateful to Grigori Olshanski for the setting of
the problem, permanent attention and fruitful discussions.
\section{Ewens-Pitman's partition structures}\label{s1}
\subsection{Notation}\label{s1.1}
In this subsection we 
give some combinatorial notation
which is used throughout the paper.

A partition is a sequence of the form
\begin{equation*}
	\la=(\la_1,\la_2,\dots,\la_\ell,0,0,\dots),\qquad \la_1\ge\la_2\ge\dots\ge\la_\ell>0,
\end{equation*}
where $\la_i\in\N$ 
and only a finite number of elements differs from zero.
Partitions are identified with Young diagrams as in \cite{Ma}.
We denote them by same letters.
The number of boxes in $\la$ is denoted by 
$|\la|=\la_1+\dots+\la_\ell$.
The number of rows in $\la$ is called the 
{\em{}length of a diagram\/}
and is denoted by $\ell(\la)$.
Also let $\varnothing$ denote the empty Young diagram.

If a diagram $\la$ 
is obtained from a diagram 
$\mu$
by adding one box, then
we write
$\mu\nearrow\la$ or, equivalently, $\la\searrow\mu$.
Denote this box 
(that distinguishes $\la$ and $\mu$)
by $\la/\mu$.

Let $\square$ be an arbitrary box.
By $\con(\square)$ denote its column number
counting from left to right.
This number does not depend on a partition
containing the box.

Fix an arbitrary nonempty diagram $\la$.
The following important properties hold.
First, for each $i=1,\dots,\ell(\la)$
there exists a unique box with $\con(\square)=\la_i$
(denoted by $\square(\la_i)$)
such that it can be removed from $\la$
and the result is again a Young diagram.
The resulting diagram is denoted by $\la-\square(\la_i)$.
Every diagram $\mu$ such that $\mu\nearrow\la$
has this form.
Second, for each $i=1,\dots,\ell(\la)$
there exists a unique box with $\con(\square)=\la_i+1$
(denoted by $\square(\la_i+1)$)
such that it can be added to $\la$
and the result is again a Young diagram.
Denote the resulting diagram by $\la+\square(\la_i+1)$.
Every diagram $\nu$ such that $\nu\searrow\la$
has either the form $\la+\square(\la_i+1)$ 
for some $i=1,\dots,\ell(\la)$ or
the form $\la+\square(1)$, where $\la+\square(1)$
is the diagram obtained from $\la$ by adding a one-box row.

For a natural $k$ by $\kr{\la}{k}$ denote
the number of rows in $\la$ of length $k$. 
It is a nonnegative integer.
For two diagrams $\mu,\la$ such that $|\la|=|\mu|+1$ we set
\begin{equation*}
	\kra{\la}{\mu}:=\left\{
	\begin{array}{ll}
		\kr{\la}{\con(\la/\mu)},&\mbox{if $\mu\nearrow\la$};\\
		0,&\mbox{otherwise}.
	\end{array}
	\right.
\end{equation*}

All Young diagrams are organized in a graded set
$\Kb:=\bigsqcup_{n=0}^{\infty}\Kb_n$,
where
$\Kb_n=\left\{ \la\colon|\la|=n \right\}$, $n\in\N$, $\Kb_0=\left\{ \varnothing \right\}$.
We introduce the structure of a graded
graph on this set. This graph has edges only between
consecutive ``floors'' $\Kb_n$ and $\Kb_{n+1}$.
If 
$\mu\in\Kb_n$ and $\la\in\Kb_{n+1}$ for some $n\ge0$, then
we draw $\kra{\la}{\mu}$ edges between $\mu$ and $\la$.
Let edges be oriented in the direction from $\Kb_n$ to $\Kb_{n+1}$.
This graph differs from the Young graph.
Namely, the latter has the same edges 
without multiplicities.

By $g(\mu,\la)$ denote
the total number of oriented paths 
from $\mu$ to $\la$ in the graph $\Kb$.
Clearly, $g(\mu,\la)$ 
vanishes unless $\mu\subset\la$
as diagrams (the sets of boxes in the plane).
Set 
$g(\la):=g(\varnothing,\la)$,
it appears that
$g(\la)=|\la|!/(\la_1!\dots\la_{\ell(\la)}!)$.

We will also need Pochhammer symbols
\begin{equation*}
	(a)_k:=a(a+1)\dots(a+k-1),\quad k=1,2,\dots,\qquad (a)_0:=1
\end{equation*}
and factorial powers
\begin{equation*}
	a^{\da k}:=a(a-1)\dots(a-k+1),\quad k=1,2,\dots,\qquad a^{\da0}:=1.
\end{equation*}

\subsection{The up/down Markov chains}\label{s1.2}
The graph structure on $\Kb$ is important, but
in this paper it is needed only to use 
the recurrent relations (\ref{f1})
for $g(\mu,\la)$ below.

One of the most important objects under consideration
are the {\em{}up/down Markov chains\/}.
They were studied in \cite{BO,Fu1,Fu2}.\footnote{The papers \cite{Fu1,Fu2}
introduced and studied 
the down/up chains,
but their difference 
from the up/down chains is minor.}
We use the formalism of \cite[\S1]{BO} to define them.

The {\em{}down transition function\/} for 
$\mu,\la\in\Kb$ such that $|\la|=|\mu|+1$ is defined as
\begin{equation*}
	p^\da(\la,\mu):=\frac{g(\mu)}{g(\la)}\kra{\la}{\mu}.
\end{equation*}
It can be easily checked that
\begin{itemize}
	\item $p^\da(\la,\mu)\ge0$ for all $\mu,\la\in\Kb$ such that $|\la|=|\mu|+1$;
	\item $p^\da(\la,\mu)$ vanishes unless $\mu\nearrow\la$;
	\item If $|\la|=n\ge1$, then $\sum_{\mu\colon|\mu|=n-1}p^\da(\la,\mu)=1$.
\end{itemize}
The object $(\Kb,p^\da)$ gives rise to {\em{}partition
structures\/}
$\{M_n\}_{n\ge0}$, where $M_n$ is a probability measure on $\Kb_n$ for every $n\ge0$
and 
\begin{equation}\label{f2}
	M_n(\mu)=\sum_{\la\colon\la\searrow\mu}M_{n+1}(\la)p^\da(\la,\mu)\quad\mbox{for all $n\in\Nn$
	and $\mu\in\Kb_n$.}
\end{equation}
Here by $M_n(\mu)$ we denote the measure of a singleton $\{ \mu \}$.

Fix a partition structure 
$\left\{ M_n \right\}$. The {\em{}up transition function\/}
for $\la,\nu\in\Kb$ such that $|\nu|=|\la|+1$ and $M_n(\la)\ne0$
is defined as
\begin{equation*}
	p^\ua(\la,\nu):=\frac{M_{n+1}(\nu)}{M_n(\la)}p^\da(\nu,\la).
\end{equation*}
This function depends on the choice of 
a partition structure. Moreover,
$\left\{ M_n \right\}$ and $p^\ua$ 
are consistent in a sense similar to
(\ref{f2}):
\begin{equation*}
	M_{n+1}(\nu)=\sum_{\la\colon\la\nearrow\nu,\; \la\in\mathop{\mathrm{supp}}(M_n)}M_n(\la)p^\ua(\la,\nu)
	\quad\mbox{for all $n\in\Nn$
	and $\nu\in\Kb_{n+1}$.}
\end{equation*}

Let $\left\{ M_n \right\}$ be a partition structure and $M_n(\la)>0$ for all $n\ge0$ and
$\la\in\Kb_n$. 
Let us define a Markov chain $T_n$ on each $\Kb_n$, $n\ge1$
with the following transition matrix:
\begin{equation*}
	T_n(\la,\widetilde \la):=\sum_{\nu\colon|\nu|=n+1}p^\ua(\la,\nu)p^\da(\nu,\widetilde\la),\qquad
	\la,\widetilde\la\in\Kb_n.
\end{equation*}
This is the composition of the up and
down transition functions, from $\Kb_n$ to $\Kb_{n+1}$ and then back to $\Kb_n$.
From the definitions above it follows that $M_n$
is a stationary distribution for $T_n$.
It can be readily shown that the matrix 
$M_n(\la)T_n(\la,\widetilde\la)$
is symmetric with respect to the substitution
$\la\leftrightarrow\widetilde\la$.
This means that the chain $T_n$ is 
reversible with respect to $M_n$.

\subsection{Ewens-Pitman's partition structures}\label{s1.3}
In the present paper we deal with a special
two-parameter family of
{\em{}Ewens-Pitman's partition structures\/}. 
It is defi\-ned as follows.

Let $\al$ and $\te$ be arbitrary parameters. 
We set
\begin{equation*}
	M_n(\la):=\frac{n!}{(\te)_n}\cdot\frac{\te(\te+\al)\dots(\te+(\ell(\la)-1)\al)}
	{\prod_{k=1}^{\infty}\kr{\la}{k}!\cdot\prod_{i=1}^{\ell(\la)}\la_i!}
	\prod_{\textstyle\genfrac{}{}{0pt}{}{\square\in\la}{\con(\square)\ge2}}(\con(\square)-1-\al).
\end{equation*}
for all $n\in\Nn$ and $\la\in\Kb_n$.

It can be proved that 
for all $\al$ and $\te$ such that $M_n$ is well-defined
$M_n$ 
satisfies
(\ref{f2}) and that for such parameters we have
\begin{equation*}
	\sum_{\la\colon|\la|=n}M_n(\la)=1\quad\mbox{for all $n\in\Nn$.}
\end{equation*}
It also can be checked that $M_n(\la)$ 
is nonnegative for all $n\ge0$ and $\la\in\Kb_n$
if and only if
\begin{itemize}
	\item ({\em{}principal series\/}) $0\le\al<1$ and $\te>-\al$; 
	\item ({\em{}degenerate series\/}) $\al=-\be<0$ and $\te=N\be$ for some $N=2,3,\dots$.
\end{itemize}
In the principal series the support of each $M_n$
is the whole $\Kb_n$. In degenerate series $M_n(\la)>0$
only for $\la$ of length $\le N$.
It follows that when the parameters are of either 
principal or degenerate series, the system $\{M_n\}_{n\ge 0}$
is a partition structure.
It is called the Ewens-Pitman's partition structure.
It can be easily shown that the up transition function
for this partition structure equals
\begin{equation*}
	p^\ua(\la,\nu)=\left\{
	\begin{array}{ll}
		\displaystyle
		\frac{\la_i-\al}{n+\te}\kr{\la}{\la_i},&
		\mbox{if $\nu=\la+\square(\la_i+1)$ for $1\le i\le\ell(\la)$;}
		\\\rule{0pt}{18pt}
		\displaystyle
		\frac{\te+\ell(\la)\al}{n+\te},&\mbox{if $\nu=\la+\square(1)$;}\\
		0,&\mbox{otherwise},
	\end{array}
	\right.
\end{equation*}
where $|\la|=n\in\Nn$.

Everywhere below except \S\ref{s5.4} we assume that
the parameters $\al$ and $\te$ are of principal series.
It was explained in \S\ref{s1.2} that 
in this case due to the positivity 
of $M_n$ on the whole $\Kb_n$
we can consider the up/down Markov chains
on $\Kb_n$, $n\ge1$.

\section{Symmetric functions and the Kingman simplex}\label{s2}
\subsection{Symmetric functions}\label{s2.1}
Let $\La$ be the (real) 
algebra of symmetric functions
in the formal variables
$y_1,y_2,\dots$.
We will need the following functions:

{\em{}Newton power sums\/}
$p_k=\sum_{i=1}^{\infty}y_i^k,$ $k=1,2,\dots$.
These elements are algebraically 
independent and generate $\La$ as a commutative unital algebra:
$\La=\mathbb{R}\left[ p_1,p_2,p_3,\dots \right]$.

{\em{}Monomial functions} $m_\la$, $\la\in\Kb$ are defined as
$\sum y_{i_1}^{\la_1}\dots y_{i_\ell(\la)}^{\la_{\ell(\la)}}$, 
where the
sum is taken over
all {\em{}distinct monomials\/},
the
indexes
$i_1,\dots,i_{\ell(\la)}$
are pairwise distinct and run from one to infinity.
We also need multiples of $m_\la$ of the form
$\m_\la:=\big(\prod_{k\ge1}\kr{\la}{k}!\big)m_\la$.
They can be viewed as similar sums
$\sum y_{i_1}^{\la_1}\dots y_{i_\ell(\la)}^{\la_{\ell(\la)}}$,
where the sum is taken over all
{\em{}collections of pairwise distinct indexes\/} $i_1,\dots,i_{\ell(\la)}$ 
from one to infinity.
It is sometimes simpler to deal with 
$\m_\la$ instead of $m_\la$.
Each of the systems $\{m_\la\}_{\la\in\Kb}$
and $\{\m_\la\}_{\la\in\Kb}$
is a basis for $\La$ as a vector space over $\mathbb{R}$.

{\em{}Factorial functions\/} $m_\la^*$ and $\m_\la^*$, $\la\in\Kb$
are obtained from $m_\la$ and $\m_\la^*$, respec\-tively, by substituting
each power of a variable $y_i^k$ by the factorial power $y_i^{\da k}$.
Hence the homogeneous component of 
$m_\la^*$ and $\m_\la^*$ of
maximal degree $|\la|$ is equal to
$m_\la$ and $\m_\la$, respectively.
It follows that each of the systems 
$\{m_\la^*\}_{\la\in\Kb}$
and $\{\m_\la^*\}_{\la\in\Kb}$ 
is also a basis for $\La$ as a vector space over $\mathbb{R}$.

Let $I:=(p_1-1)\La$ be the principal
ideal in $\La$ generated by $p_1-1$.
Set $\Lao:=\La/I$.
To every element $f\in\La$ corresponds an image in $\Lao$ denoted by
$f^\circ$. In particular, $p_1^\circ=1$
and $\Lao$ is freely generated (as a commutative unital algebra)
by the elements $p_k^\circ$, $k=2,3,\dots$.
Moreover,
\begin{equation*}
	\La=\mathbb{R}\left[ p_1,p_2,p_3\dots \right]=
	I\oplus\mathbb{R}\left[ p_2,p_3,\dots \right].
\end{equation*}
It follows that $\Lao\cong\mathbb{R}\left[ p_2,p_3,\dots \right]$.
It can be easily checked that
the basis for the latter algebra over $\mathbb{R}$ is
$\{m_\la\}_{\la\in\Kb,\; \kr{\la}1=0}$.
Therefore,
the basis for $\Lao$ over $\mathbb{R}$ is
$\{m_\la^\circ\}_{\la\in\Kb,\; \kr{\la}1=0}$.

\subsection{The Kingman simplex and moment coordinates}\label{s2.2}
In the Introduction we described the {\em{}Kingman simplex\/}
\begin{equation*}
	\sk=\biggl\{x=(x_1,x_2,\dots)\in[0,1]^\infty\colon x_1\ge x_2\ge\dots\ge0,\
	{\displaystyle \sum\limits_{i=1}^\infty x_i\le 1}\biggr\}.
\end{equation*}
It is a compact, metrizable
and separable space
in the topology of coordinatewise convergence
(we use this topology throughout the paper).
The simplex $\sk$ contains a distinguished dense subspace
$\nabla_\infty:=\{x=(x_1,x_2,\dots)\in\sk\colon\sum_{i=1}^{\infty}x_i=1\}$.
The use of the symbols $\nabla_\infty$ and $\sk$ follows the work \cite{EK1}.
By $C(\sk)$ denote the algebra of real continuous
functions on $\sk$ with pointwise operations and the supremum norm.

To every point $x\in\sk$ we assign a probability measure
\begin{equation*}
	\nu_x:=\sum_{i=1}^{\infty}x_i\delta_{x_i}+\gamma(x)\delta_0
\end{equation*}
on the segment
$\left[ 0,1 \right]$, where $\delta_s$ 
is the Dirac measure at a point $s$ and $\gamma(x):=1-\sum_{i=1}^{\infty}x_i$. 
By $q_k(x)$ 
denote the $k$th moment of $\nu_x$:
\begin{equation*}
	q_k(x):=\int_0^1u^k\nu_x(du)=\sum_{i=1}^{\infty}x_i^{k+1},\qquad k=1,2,\dots.
\end{equation*}
These functions are continuous on $\sk$ because
$x_i\le i^{-1}$ for every
$i=1,2,\dots$.
It is worth noting that the function $\gamma(x)$ is not continuous on $\sk$.
The functions $q_k(x)$, $k\ge1$ also separate points of $\sk$ 
(because the measure on $\left[ 0,1 \right]$ 
in uniquely determined by its moments) and are algebraically independent.
Following \cite{BO} we call
$q_1(x),q_2(x),\dots$ 
the {\em{}moment coordinates\/}
of a point $x\in\sk$.

Let $\Fc=\mathbb{R}\left[ q_1,q_2,\dots \right]$
be the commutative unital algebra generated by the moment coordinates.
By the Stone-Weierstrass theorem,
$\Fc$ is a dense subalgebra of $C(\sk)$ .

\subsection{Symmetric functions on the simplex. The Kingman theorem}\label{s2.3}
The correspondence $p_2^\circ\to q_1(x),
p_3^{\circ}\to q_2(x),\dots$
establishes an isomorphism of the algebra $\Lao$ described in \S\ref{s2.1}
and the algebra $\Fc$ from \S\ref{s2.2}.
Hence
to each element $f^\circ\in\Lao$ corresponds
a continuous function on $\sk$. Denote this function by $f^\circ(x)$.
In particular,
$p_k^\circ(x)=\sum_{i=1}^{\infty}x_i^k$, $k=2,3,\dots$ and $p_1^\circ(x)\equiv1$.

Now we recall the Kingman theorem about partition structures.
Consider for $n=1,2,\dots$ the following embeddings
\begin{equation*}
	\iota_n\colon\Kb_n\hookrightarrow\sk,\qquad
	\iota_n\colon(\la_1,\dots,\la_\ell)\mapsto
	\left( \frac{\la_1}n,\dots,\frac{\la_\ell}n,0,0,\dots \right)\in \sk,
\end{equation*}
where $\la=(\la_1,\dots,\la_\ell)\in\Kb_n$.
The next remark will be useful in \S\ref{s4}.
\begin{note}\label{p2.1}\rm{}
        The sets $\iota_n(\Kb_n)$ approximate the space $\sk$
	in the sense that any open subset of $\sk$ 
	has a nonempty intersection with $\iota_n(\Kb_n)$ 
	for all $n$ large enough.
\end{note}

\begin{thm}[\cite{Ki3,KOO}]\label{p2.2}
	For every partition structure 
	$\left\{ M_n \right\}$ on $\Kb$ there exists a
	Borel probability measure $P$ 
	on $\sk$ (called the {\em{}boundary measure\/} of a partition structure)
	such that
	\begin{equation*}
		P=\lim_{n\to\infty}\iota_n(M_n).
	\end{equation*}
	Conversely,
	any partition structure can be
	reconstructed from its boundary measure as follows:
	\begin{equation*}
		M_n(\la)=g(\la)\int_{\sk}m_\la^\circ(x)P(dx)\quad\mbox{for all $\la\in\Kb_n$.}
	\end{equation*}
	Hence the partition structures on $\Kb$ and Borel probability measures
	on $\sk$ are in one-to-one correspondence.
\end{thm}

To the Ewens-Pitman's partition structure (\S\ref{s1.3})
with parameters $\al$ and $\te$ of principal or degenerate series
corresponds the  
well-known
{\em{}Poisson-Dirichlet measure\/}
$\pd{\al}{\te}$. 
We have\footnote{The functions $\m_\la$ are used instead of $m_\la$ to simplify the notation.}
\begin{equation*}
	\int_{\sk}\m_\la^\circ(x)\,\pd\al\te(dx)=\frac{(-\te/\al)^{\da\ell(\la)}}{(\te)_{|\la|}}
	\prod_{i=1}^{\ell(\la)}(-\al)_{\la_i}\quad\mbox{for all $\la\in\Kb$}.
\end{equation*}
When $\al=0$, the value of the RHS is obtained by continuity.

It is known (see, e.g., \cite{PPY,Pi-Rand})
that the measure $\pd{\al}{\te}$ is concentrated on
the dense subset $\nabla_\infty\subset\sk$ 
(the definition of $\nabla_\infty$ see in \S\ref{s2.2}).

\section{The up/down Markov chains transition operators' action on symmetric functions}\label{s3}
Here we assume that the parameters $\al$ and $\te$ are of
principal series.
Let us introduce some extra notation first.

For every set $\mathcal{X}$
denote by $\mathrm{Fun}(\mathcal{X})$
the algebra of real functions on $\mathcal{X}$
with pointwise operations.

Consider an embedding of the algebra 
of symmetric functions $\La$ 
into the algebra $\mathrm{Fun}(\Kb)$
defined on the generators $p_k$, $k=1,2,\dots$ as follows:
$p_k\to p_k(\la):=\sum_{i=1}^{\ell(\la)}\la_i^k$.
Thus, to every element $f\in\La$ corresponds
a function from $\mathrm{Fun}(\Kb)$. Denote this function by $f(\la)$.

Let $f\in\La$. By $f_n$ denote the restriction of 
the function $f(\cdot)\in\mathrm{Fun}(\Kb)$ to $\Kb_n\subset\Kb$.
It can be easily checked that the subalgebra $\La\subset\mathrm{Fun}(\Kb)$
separates points. Therefore, the functions of the form $f_n$,
with $f\in\La$, exhaust the (finite-dimensional) 
space $\mathrm{Fun}(\Kb_n)$.

The aim of this section is to prove the 
following proposition which is needed for Lemma \ref{p4.1} below.
\begin{prop}\label{p3.1}
	Consider the transition operator of
	the $n$th up/down Markov chain
	$T_n\colon\mathrm{Fun}(\Kb_n)\to\mathrm{Fun}(\Kb_n)$
	which corresponds to the
	two-parameter Ewens-Pitman's partition structure.
	Its action on the functions
	$\m_\mu^*$, $\mu\in\Kb$
	looks as follows:
	\begin{equation*}
		\left.
		\begin{array}{l}\displaystyle
			(T_n-{\bf1})(\m_\mu^*)_n=-\frac{k(k-1+\te)}{(n+1)(n+\te)}(\m_\mu^*)_n
			\\\displaystyle\qquad+
			\frac{n+1-k}{(n+1)(n+\te)}
			\sum_{\textstyle\genfrac{}{}{0pt}{}{i=1}{\mu_i\ge2}}^{\ell(\mu)}
			\mu_i\big(\mu_i-1-\al\big)
			\big(\m_{\mu-\square(\mu_i)}^*\big)_n\\\displaystyle\qquad+
			\frac{n+1-k}{(n+1)(n+\te)}
			\krbig{\mu}{1}\big(\te+\al(\ell(\mu)-1)\big)\big(\m^*_{\mu-\square(1)}\big)_n,
		\end{array}
		\right.
	\end{equation*}
	where ${\bf1}$ 
	denotes the identity operator and $k=|\mu|$.
\end{prop}
\begin{note}\label{p3.2}\rm{}
The Proposition states that
$(T_n-{\bf1})(\m_\mu^*)_n$
is a linear combination of the functions
$(\m_\mu^*)_n$ 
(the first summand in the RHS)
and the functions of the form
$(\m_\varkappa^*)_n$
for all
$\varkappa\colon\varkappa\nearrow\mu$.
The sum over $i$
in the RHS deals with 
$\varkappa$ of the same length as
$\mu$. 
Note that this sum contains similar terms if
$\mu$ 
has equal rows of length $\ge2$.
The last summand is nonzero if $\kr{\mu}{1}>0$
and in this case corresponds to the diagram
$\varkappa=\mu-\square(1)$
obtained from $\mu$
by deleting a one-box row.
\end{note}

Let us give two formulas
concerning factorial functions and the 
numbers of paths $g(\mu,\la)$ defined in \S\ref{s1.1}.
We use these formulas below to proof the
Proposition \ref{p3.1}.
Let
$|\la|=n\ge m=|\mu|$, $\la,\mu\in\Kb$.
Then
\begin{equation}\label{f1.1}
	\frac{g(\mu,\la)}{g(\la)}=\frac{(m_\mu^*)_n(\la_1,\dots,\la_{\ell(\la)})}
	{n(n-1)\dots(n-m+1)},
\end{equation}
and
\begin{equation}\label{f1}
	g(\mu,\la)=\sum_{\varkappa\colon\varkappa\nearrow\la}\kra{\la}{\varkappa}
	g(\mu,\varkappa).
\end{equation}
The formula (\ref{f1.1}) can be checked directly, and
(\ref{f1}) 
follows from the definitions of
$g(\mu,\la)$ and $\kappa(\varkappa,\la)$,
see \S\ref{s1.1}.

In order to prove Proposition
\ref{p3.1} let us write the operator
$T_n$ 
as a composition of ``down''
$D_{n+1,n}\colon \mathrm{Fun}(\Kb_n)\to\mathrm{Fun}(\Kb_{n+1})$
and ``up''
$U_{n,n+1}\colon\mathrm{Fun}(\Kb_{n+1})\to\mathrm{Fun}(\Kb_{n})$
operators acting on functions.

The operator $D_{n+1,n}$ 
is constructed using the down transition probabilities
and does not depend on the parameters $\al$ and $\te$.
The operator
$U_{n,n+1}$ is constructed using the up transition probabilities
and depends on the parameters. Namely,
\begin{equation}\label{f2.9}
	\left.
	\begin{array}{rcll}
		(D_{n+1,n}f_n)(\la)&:=& \displaystyle\sum_{\mu\colon\mu\nearrow\la}p^\da(\la,\mu)f_n(\mu),&\qquad\la\in\Kb_{n+1}\\
		\rule{0pt}{16pt}
		(U_{n,n+1}f_{n+1})(\mu)&:=& \displaystyle\sum_{\la\colon\la\searrow\mu}p^\ua(\mu,\la)f_{n+1}(\la),&\qquad\mu\in\Kb_n.
	\end{array}
	\right.
\end{equation}
These operators are adjoint to the corresponding operators
acting on measures.
The latter act in agreement with their names, e.g.,
$D_{n+1,n}^*\colon
\mathcal{M}(\Kb_{n+1})\to
\mathcal{M}(\Kb_{n})$,
where $\mathcal{M}(\mathcal{X})$ is the space of measures on $\mathcal{X}$.

It clearly follows from the definition of the $n$th up/down Markov chain
(\S\ref{s1.2})
that
$T_n=U_{n,n+1}\circ D_{n+1,n}\colon\mathrm{Fun}(\Kb_n)\to\mathrm{Fun}(\Kb_n)$, $n\in\N$.
Proposition \ref{p3.1} 
directly follows from
\begin{lemma}
	{\rm{}(1)\/} 
	There exists a unique operator $\widetilde D\colon\La\to\La$ such that
	\begin{equation*}
		D_{n+1,n}f_n=\frac1{n+1}(\widetilde Df)_{n+1}
	\end{equation*}
	for all $n\in\Nn$ and $f\in\La$.
	In the basis $\{ \m_\mu^* \}_{\mu\in\Kb}$ for the algebra $\La$ 
	this operator has the form
	\begin{equation}\label{f3}
		\widetilde D\m_\mu^*=(p_1-|\mu|)\m_\mu^*.
	\end{equation}

	{\rm{}(2)\/} There exists a unique operator $\widetilde U\colon\La\to\La$ such that
	\begin{equation*}
		U_{n,n+1}g_{n+1}=\frac1{n+\theta}(\widetilde Ug)_{n}
	\end{equation*}
	for all $n\in\Nn$ and $g\in\La$.
	In the basis $\{\m_\mu^*\}_{\mu\in\Kb}$ for the algebra $\La$ 
	this operator has the form
	\begin{equation}\label{f4}
		\left.
		\begin{array}{r}
			\displaystyle
			\widetilde U\m_\mu^*=
			\big(p_1+\theta+|\mu|\big)\m_\mu^*
			+\sum_{\textstyle\genfrac{}{}{0pt}{}{i=1}{\mu_i\ge2}}^{\ell(\mu)}
			\mu_i\big( \mu_i-1-\alpha \big)\m_{\mu-\square(\mu_i)}^*\phantom{.}
			\\\displaystyle\rule{0pt}{12pt}
			+\krbig{\mu}{1}\big(\theta+\alpha(\ell(\mu)-1)\big)\m_{\mu-\square(1)}^*.
		\end{array}
		\right.
	\end{equation}
\end{lemma}
Note that in the RHS of (\ref{f4})
there is a sum similar to the one explained in Remark
\ref{p3.2}.
\begin{proof} (1) Let us show that the operator
	$\widetilde D$ defined by (\ref{f3})
	is the desired one. 
	We use the connection (\ref{f1.1}) 
	of factorial functions with the numbers of paths
	$g(\mu,\la)$ and the recurrent relations (\ref{f1}) for the latter.

	Let $\mu\in\Kb$, $\nu\in\Kb_{n+1}$. We suppose that 
	$n\ge|\mu|$, because otherwise  $(m_\mu^*)_n=0$.
	We have
	\begin{equation*}
		\left.
		\begin{array}{l}\displaystyle
			(D_{n+1,n}(m_\mu^*)_n)(\nu)=
			\sum_{\la\colon\la\nearrow\nu}p^\da(\nu,\la)(m_\mu^*)_n(\la)\\\displaystyle\quad=
			\sum_{\la\colon\la\nearrow\nu}\frac{g(\la)}{g(\nu)}
			\kra{\nu}{\la}(m_\mu^*)_n(\la)=
			\sum_{\la\colon\la\nearrow\nu}\frac{g(\la)}{g(\nu)}
			\kra{\nu}{\la}n^{\da|\mu|}\frac{g(\mu,\la)}{g(\la)}\\\displaystyle\quad=
			\frac{n^{\da|\mu|}}{g(\nu)}\sum_{\la\colon\la\nearrow\nu}
			\kra{\nu}{\la}g(\mu,\la)=n^{\da|\mu|}\frac{g(\mu,\nu)}{g(\nu)}=
			\frac{n+1-|\mu|}{n+1}(m_\mu^*)_{n+1}(\nu).
		\end{array}
		\right.
	\end{equation*}
	Thus,
	\begin{equation*}
		\frac1{n+1}\widetilde Dm_\mu^*(\nu)=
		\frac1{n+1}\big((p_1-|\mu|)m_\mu^*\big)(\nu),
	\end{equation*}
	because $p_1(\nu)=|\nu|=n+1$.
	It follows that in the basis $\{ m_\mu^*\}_{\mu\in\Kb}$ 
	the operator $\widetilde D$ has the form
	$\widetilde Dm_\mu^*=(p_1-|\mu|)m_\mu^*$. 
	If we multiply both sides by $\prod_{k\ge1}\kr{\mu}{k}!$, then
	we get the desired expression in the basis $\{\m_\mu^*\}_{\mu\in\Kb}$.
	The uniqueness of the operator $\widetilde D$ follows from the fact that $\La$
	is embedded into $\mathrm{Fun}(\Kb)$.

	(2) Fix $\mu\in\Kb$ and denote $l:=\ell(\mu)$. Let $\nu\in\Kb_n$.
	We have an explicit expression
	\begin{equation*}
		(\m_\mu^*)_n(\nu)=\sum\nu_{j_1}^{\da\mu_1}\dots\nu_{j_l}^{\da\mu_l},
	\end{equation*}
	where the sum is taken over all pairwise distinct $j_1,\dots,j_l$ from $1$ to $\ell(\nu)$.
	Thus, we can verify (\ref{f4}) directly.
		
	From the definition of the operator $U_{n,n+1}$ (\ref{f2.9})
	and the formula for the 
	up transition probabilities of the 
	Ewens-Pitman's partition structure
	(\S\ref{s1.3})
	it follows that
	\begin{equation}\label{f5}
		\left.
		\begin{array}{rcr}
			\displaystyle
			\big(U_{n,n+1}(\m_\mu^*)_{n+1}\big)(\nu)&=& \displaystyle
			\frac1{n+\te}\Bigg[
			\sum_{i=1}^{\ell(\nu)}
			(\nu_i-\al)(\m_\mu^*)_{n+1}\big(\nu+\square(\nu_i+1)\big)\phantom{.}\\&&\displaystyle
			\phantom{\frac1{n+\te}\Bigg[}
			+\big(\te+\ell(\nu)\al\big)(\m^*_\mu)_{n+1}\big(\nu+\square(1)\big)
			\Bigg].
		\end{array}
		\right.
	\end{equation}
        Let us transform the RHS and obtain (\ref{f4}). 

	{\bf1.} First, note two simple properties:
	\begin{equation*}	
		\mbox{\bf1a. } (a+1)^{\da b}-a^{\da b}=b\cdot
		a^{\da(b-1)}; \qquad\mbox{\bf1b. } a\cdot a^{\da
		b}=a^{\da (b+1)}+b\cdot a^{\da b}.
	\end{equation*}

	{\bf2.} Let us deal with the standalone summand in the RHS of (\ref{f5}).
	We show that
	\begin{equation}\label{f6}
		(\m_\mu^*)_{n+1}\big(\nu+\square(1)\big)=
		(\m_\mu^*)_n(\nu)+\kr{\mu}{1}\big(\m_{\mu-\square(1)}^*\big)_n(\nu).
	\end{equation}
	It is clear that $(\m_\mu^*)_{n+1}\big(\nu+\square(1)\big)=
	\sum(\nu+\square(1))_{j_1}^{\da\mu_1}\dots(\nu+\square(1))_{j_l}^{\da\mu_l}$,
	where the sum is taken over all pairwise distinct
	$j_1,\dots,j_l$ from $1$ to $\ell(\nu)+1$.
	To every combination of indexes
	$j_1,\dots,j_l$ we assign a number
	$c\in\left\{ 1,\dots,l \right\}$ defined by the condition $j_c=\ell(\nu)+1$.
	It follows that $(\nu+\square(1))_{j_c}^{\da\mu_c}=1^{\da\mu_c}$,
	and for other $j_k$ we have $(\nu+\square(1))_{j_k}=\nu_{j_k}$.
	Let us combine the summands of the form 
	$(\nu+\square(1))_{j_1}^{\da\mu_1}\dots(\nu+\square(1))_{j_l}^{\da\mu_l}$
	with equal $c$, i.e., let us write
	\begin{equation*}
		\begin{array}{l}\displaystyle
			(\m_\mu^*)_{n+1}\big(\nu+\square(1)\big)=
			(\m_\mu^*)_n(\nu)
			\\\displaystyle\phantom{(\m_\mu^*)_{n+1}\big(\nu+\square(1)\big)=}\rule{0pt}{22pt}
			+\sum_{c=1}^{l}1^{\da\mu_c}
			\Big(
			\sum_{\mbox{\scriptsize{}w/o $j_c$}}
			\nu_{j_1}^{\da\mu_1}\dots
			\nu_{j_{c-1}}^{\da\mu_{c-1}}
			\nu_{j_{c+1}}^{\da\mu_{c+1}}
			\dots
			\nu_{j_l}^{\da\mu_{j_l}}
			\Big).
		\end{array}
	\end{equation*}
	The last sum
	$\sum_{\mbox{\scriptsize{}w/o $j_c$}}$
	means the sum taken over all pairwise distinct indexes $j_1,\dots,j_{c-1},
	j_{c+1},\dots,j_l$ from $1$ to $\ell(\nu)$. 
	We denote this sum by $S_c$ to simplify the notation.
	Note that $1^{\da b}=0$ if $b\ge2$, and $1^{\da1}=1$. Therefore,
	\begin{equation*}
		\sum_{c=1}^{l}1^{\da\mu_c}S_c=\kr{\mu}{1}\big(\m_{\mu-\square(1)}^*\big)_n(\nu),
	\end{equation*}
	and this implies (\ref{f6}).

	{\bf3.} 
	Now we deal with the sum over $i$ in the RHS of (\ref{f5}).

	First, fix an arbitrary 
	$i$ from $1$ to $\ell(\nu)$. We have
	\begin{equation*}
		\left.
		\begin{array}{lcll}
			(\m_\mu^*)_n(\nu)&=& \displaystyle
			\sum_{c=1}^{l}\nu_i^{\da\mu_c}S_c+
			\sum\limits_{j_1,\dots,j_l\ne i}
			\nu_{j_1}^{\da\mu_1}\dots\nu_{j_l}^{\da\mu_l};\\
			(\m_\mu^*)_{n+1}\big(\nu+\square(\nu_i+1)\big)&=& \displaystyle
			\sum_{c=1}^{l}(\nu_i+1)^{\da\mu_c}S_c+
			\sum\limits_{j_1,\dots,j_l\ne i}
			\nu_{j_1}^{\da\mu_1}\dots\nu_{j_l}^{\da\mu_l}.
		\end{array}
		\right.
	\end{equation*}
	Using {\bf1a}, we get
	\begin{equation*}
		(\m_\mu^*)_{n+1}\big(\nu+\square(\nu_i+1)\big)-
		(\m_\mu^*)_n(\nu)=\sum_{c=1}^{l}\mu_c\nu_i^{\da(\mu_c-1)}S_c.
	\end{equation*}
	Thus, the sum over $i$ in the RHS of
	(\ref{f5}) becomes
	\begin{equation*}
		\left.
		\begin{array}{l}\displaystyle
			\sum_{i=1}^{l}(\nu_i-\al)(\m_\mu^*)_{n+1}\big(\nu+\square(\nu_i+1)\big)\\\displaystyle\qquad=
			\big(|\nu|-\ell(\nu)\al\big)(\m_\mu^*)_n(\nu)+
			\sum_{i=1}^{\ell(\nu)}\sum_{c=1}^{l}(\nu_i-\al)\mu_c\nu_i^{\da(\mu_c-1)}S_c.
		\end{array}
		\right.
	\end{equation*}
	Now fix an arbitrary $c$ from $1$ to $l$. Using {\bf1b}, we get
	\begin{equation*}
		\left.
		\begin{array}{l}\displaystyle
			\sum_{i=1}^{\ell(\nu)}(\nu_i-\al)\mu_c\nu_i^{\da(\mu_c-1)}S_c\\\displaystyle\qquad=
			S_c\Bigg[
			\mu_c\sum_{i=1}^{\ell(\nu)}\nu_i^{\da\mu_c}+
			\mu_c(\mu_c-1)\sum_{i=1}^{\ell(\nu)}\nu_i^{\da(\mu_c-1)}
			-\al\mu_c\sum_{i=1}^{\ell(\nu)}\nu_i^{\da(\mu_c-1)}
			\Bigg]\\\rule{0pt}{22pt}\displaystyle\qquad=\left\{
			\begin{array}{ll}\displaystyle
				\mu_c(\m_\mu^*)_n(\nu)+\mu_c(\mu_c-1-\al)
				\big(\m_{\mu-\square(\mu_c)}\big)_n(\nu),&
				\mbox{if }\mu_c\ge2;\\\rule{0pt}{14pt}
				(\m_\mu^*)_n(\nu)-\al\big(\ell(\nu)-(l-1)\big)
				\big(\m_{\mu-\square(1)}^*\big)_n(\nu),&\mbox{if }\mu_c=1.
			\end{array}
			\right.
		\end{array}
		\right.
	\end{equation*}

	{\bf4.} 
	If we put together {\bf2} and {\bf3} and recall that
	$|\nu|=n$ and $\ell(\mu)=l$, we get
	\begin{equation*}
		\left.
		\begin{array}{l}\displaystyle
			\big(U_{n,n+1}(\m_\mu^*)_{n+1}\big)(\nu)=
			(\m_\mu^*)_n(\nu)+\frac1{p_1(\nu)+\te}\Big[
			|\mu|(\m_\mu^*)_n(\nu)
			\\\displaystyle\qquad+
			\sum_{\textstyle\genfrac{}{}{0pt}{}{c=1}{\mu_c\ge2}}^{\ell(\mu)}
			\mu_c(\mu_c-1-\al)\big(\m^*_{\mu-\square(\mu_c)}\big)_n(\nu)\\\displaystyle\qquad+
			\kr{\mu}{1}\big(\te+\al(\ell(\mu)-1)\big)
			\big(\m^*_{\mu-\square(1)}\big)_n(\nu)
			\Big].
		\end{array}
		\right.
	\end{equation*}
	This coincides with the desired expression (\ref{f4}).
	The uniqueness of the operator $\widetilde U$ again follows
	from the fact that $\La\hookrightarrow\mathrm{Fun}(\Kb)$ is embedding.
\end{proof}

\section{Convergence of the up/down Markov chains}\label{s4}
In this section we consider 
the limit, as $n\to\infty$, of the Markov chains $T_n$
on $\Kb_n$. We assume that the parameters $\al$ and $\te$ are of principal series.
The limit is a continuous time Markov process on $\sk$ denoted by $\x{\al}{\te}(t)$.

The argument in this section 
uses the results of the book \cite{EK2}.
The application of these results is
based on the algebraic calculations
of \S\ref{s3} and is similar to \cite[\S1]{BO}.
Due to this similarity the proofs 
repeating those from \cite{BO} are not given.

\subsection{The construction of the process in the simplex}\label{s4.1}

First, let us introduce essential notation.

Each operator $T_n$ acts in a finite-dimensional space
of functions
$\mathrm{Fun}(\Kb_n)$, $n=1,2,\dots$.
These spaces can be viewed as Banach spaces with the 
supremum norm
$\|\cdot\|_n$.
Recall the embeddings 
$\iota_n\colon\Kb_n\hookrightarrow\sk$
introduced in \S\ref{s2.3}.
Let 
$\pi_n$ 
denote the corresponding projections of function spaces:
\begin{equation*}
	\big(\pi_n(f)\big)(\la):=f(\iota_n(\la)),\qquad\la\in\Kb_n,\qquad f\in C(\sk).
\end{equation*}
In \S\ref{s2} we introduced
a dense subalgebra $\Fc\subset C(\sk)$.
This algebra $\Fc$ 
admits an ascending filtration by finite-dimensional subspaces
\begin{equation*}
	\Fc^0\subset\Fc^1\subset\Fc^2\dots\subset\Fc,\qquad \bigcup_{m=0}^{\infty}\Fc^m=\Fc.
\end{equation*}
Recall that
$\Fc\cong\Lao=\La/(p_1-1)\La$.
We define the filtration
$\left\{ \Fc^m \right\}$ 
of $\Lao$ as the image of the filtration of the algebra $\La$ 
by degrees of polynomials
(in the formal variables $y_i$, see \S\ref{s2.1}).

From Remark
\ref{p2.1} 
and the fact
$\dim\Fc^m<\infty$ it follows that
each projection 
$\pi_n$ 
is one-to-one 
on $\Fc^m$
(for fixed $m$)
for all $n$ large enough.

Now we prove that the generators of the up/down Markov chains
converge to some operator $A$ in the space $C(\sk)$.
\begin{lemma}\label{p4.1}
	Each $\pi_n(\Fc^m)$ 
	is invariant under the operator $T_n$
	(for fixed $m$) for all $n$ large enough.
	Under the identification $\Fc^m=\pi_n(\Fc^m)$ there exists a limit
	\begin{equation}\label{f9}
		\lim_{n\to\infty}n^2(T_n-{\bf1})f=Af\quad\mbox{for all $f\in\Fc$}
	\end{equation}
	in any finite-dimensional space $\Fc^m$.\footnote{Here by ${\bf1}$ we denote the identity operator.}
	In this way we obtain an operator $A\colon\Fc\to\Fc$ with the property $A\Fc^m\subset\Fc^m$, $m\ge0$.
	The action of this operator
	on the functions 
	$\m_\mu^\circ\in\Lao\cong\Fc$, $\mu\in\Kb$ (they were defined in \S\ref{s2.1})
	has the following explicit form:
	\begin{equation}\label{f10}
		\begin{array}{r}\displaystyle
			A\m_\mu^\circ=-|\mu|(|\mu|-1+\te)\m_\mu^\circ+
			\sum_{\textstyle\genfrac{}{}{0pt}{}{c=1}{\mu_c\ge2}}^{\ell(\mu)}
			\mu_c(\mu_c-1-\al)\m^\circ_{\mu-\square(\mu_c)}\phantom{.}\\\displaystyle+
			\kr{\mu}{1}\big(\te+\al(\ell(\mu)-1)\big)\m_{\mu-\square(1)}^{\circ}.
		\end{array}
	\end{equation}
\end{lemma}
Note that $A$ sends the function $\m_\varnothing^\circ\equiv1$ to $0$,
as it should be.
\begin{proof}
	It is clear that for all $n$ large enough the space
	$\pi_n(\Fc^m)$ 
	is the linear span of the functions
	$(\m_\mu^*)_n$ with $|\mu|\le m$.
	Therefore, the invariance of
	$\pi_n(\Fc^m)$ under $T_n$
	follows from Proposition \ref{p3.1}, 
	because it is clear 
	that 
	the operator $T_n$ does not increase the degree of each function 
	$(\m_\mu^*)_n\in\mathrm{Fun}(\Kb_n)$
	(this degree is equal to $|\mu|$).

	Thus, identifying 
	$\pi_n(\Fc^m)$ with $\Fc^m$ (which makes sense for fixed $m$ and $n$
	large enough), we may say that each
	$\Fc^m$ is invariant under $T_n$. 
	Now we can prove 
	(\ref{f9})
	and (\ref{f10}) together.

	Consider the map $\La\to\mathrm{Fun}(\Kb_n)$ defined as
	\begin{equation*}
		f_{\left[ n \right]}(\la):=\pi_n(f^\circ)(\la)=f^{\circ}
		\left( \frac{\la_1}n,\frac{\la_2}n,\dots,\frac{\la_{\ell(\la)}}n,0,0,\dots \right),\quad
		\la\in\Kb_n,\ f\in\La.
	\end{equation*}
	This is the restriction of an element $f^{\circ}\in C(\sk)$ to
	$\iota_n(\Kb_n)\subset \sk$.
	Let $G_s\colon\La\to\La$ for all $s>0$ be an automorphism of the algebra defined 
	on its basis as
	\begin{equation*}
		G_s\m_\mu=s^{|\mu|}\m_\mu,\qquad \mu\in\Kb.
	\end{equation*}
	On the homogeneous component of degree
	$k$, $k=0,1,2,\dots$ this automorphism
	reduces to multiplication by the number
	$s^k$.\footnote{The homogeneous component of degree $k$ in the algebra $\La$
	consists of all homogeneous symmetric functions
	of degree $k$. This set is the linear span
	of the elements
	$\m_\mu$, $|\mu|=k$.}

	Therefore, we have the following expression for
	$f_n$, where $f\in\La$:
	\begin{equation}\label{f11}
		f_n(\la)=(G_nf)_{\left[ n \right]}(\la),\qquad \la\in\Kb_n.
	\end{equation}
	Indeed, since $\iota_n(\Kb_n)\subset\nabla_\infty$, it follows that
	for all $g\in\La$ and $\la\in\Kb_n$ the value
	$g^\circ(\iota_n(\la))=g_{\left[ n \right]}(\la)$
	is equal to the formal 
	evaluation of the element $g\in\La$
	on the
	coordinates
	of the point $\iota_n(\la)=(\la_1/n,\dots,\la_{\ell(\la)}/n,0,0,\dots)\in\nabla_\infty$
	(see Remark \ref{p5.1.5} below).
	In \S\ref{s2}
	we pointed out that the homogeneous component of maximal degree of 
	$\m_\mu^*$ is $\m_\mu$.
	Thus, for all $\mu\in\Kb$ we have
	\begin{equation*}
		\lim_{n\to\infty}n^{-|\mu|}G_n\m_\mu^*=\m_\mu.
	\end{equation*}
	Now, to conclude the proof of the Lemma, 
	express each function of the form $f_n(\la)$
	in the formula for
	$(T_n-{\bf1})(\m_\mu^*)_n$
	from Proposition \ref{p3.1}
	in terms of (\ref{f11}).
	It remains to 
	multiply the result by $n^{-|\mu|}$ and take the limit as $n\to\infty$.
\end{proof}
Note that the functions
$\m_\mu^\circ$, $\mu\in\Kb$ 
are not linearly independent.
However, from the proof of the Lemma it follows that the operator $A$ is well-defined by
(\ref{f10}).

Now proceed to the convergence of discrete semigroups
$\{ {\bf1},T_n,T_n^2,\dots \}$
to a {\em{}conservative Markov semigroup\/} $\{T(t)\}_{t\ge0}$
in the Banach space $C(\sk)$
constructed using the operator $A$.
Recall that a conservative Markov semigroup 
is a strongly continuous semigroup of contraction operators 
in $C(\sk)$
preserving positive functions and the constant
$1$.

We use the definition from 
\cite[Chapter 1, Section 6]{EK2}
to introduce strict sense
to the concept of
convergence of discrete semigroups to a continuous one:
\begin{df}\label{p4.2}\rm{}
	We say that a sequence of functions
	$\{f_n\in \mathrm{Fun}(\Kb_n)\}$
	converges to a function $f\in C(\sk)$ if
	$\|f_n-\pi_n(f)\|_n\to0$ as $n\to\infty$.
	Here $\|\cdot\|_n$ is the supremum norm of the space
	$\mathrm{Fun}(\Kb_n)$.
	In this case we write $f_n\to f$.
\end{df}
To establish a convergence of the semigroups
$\{ {\bf1},T_n,T_n^2,\dots\}$
to $\{T(t)\}_{t\ge0}$ 
one must perform a 
natural scaling of time: one step of the $n$th Markov chain 
corresponds to a small time interval of order $n^{-2}$.
This convergence is established in
\begin{prop}\label{p4.3}
	{\rm{}(1)\/} The operator $A\colon\Fc\to\Fc$ defined in Lemma $\ref{p4.1}$
	is closable in the space $C(\sk)$;
	
	{\rm{}(2)\/} The closure $\overline A$ of the operator $A$ generates
	a conservative Markov semigroup $\{T(t)\}_{t\ge0}$ in $C(\sk)$;

	{\rm{}(3)\/} Discrete semigroups $\{{\bf1},T_n,T_n^2,\dots\}$ converge, as $n\to\infty$,
	to the semigroup $\{T(t)\}_{t\ge0}$ in the following sense:
	\begin{equation*}
		T_n^{[ n^2t ]}\pi_n(f)\to T(t)f\quad\mbox{for all $f\in C(\sk)$}
	\end{equation*}
	(the limit in understood according to Definition \ref{p4.2})
	for all $t\ge0$ uniformly on bounded intervals.
\end{prop}
This Proposition can be proved similarly to \cite[Proposition 1.4]{BO} and
follows from the convergence of the generators
(Lemma \ref{p4.1}).
The proof uses general statements
\cite[Chapter 1, Theorem 6.5]{EK2} and
\cite[Chapter 1, Lemma 2.11]{EK2}.

Next, it directly follows from \cite[Chapter 4, Theorem 2.7]{EK2}
that $\{T(t)\}_{t\ge0}$ 
is a semigroup corresponding to a strong Markov process 
with c\`adl\`ag sample paths which can start from 
any point and any probability distribution.

We will call the operator $A$ the {\em{}pre-generator\/} of both the semigroup $\{T(t)\}_{t\ge0}$
and the process $\x{\al}{\te}(t)$. 
Thus, for every pair of parameters $\al$, $\te$
of principal series we have constructed a Markov process
$\x{\al}{\te}(t)$ in $\sk$
which can start from 
any point and any probability distribution.

\subsection{Some properties of the process in the simplex}\label{s4.2}
In this subsection we formulate and comment the 
properties of the process 
$\x{\al}{\te}(t)$
which are similar to the ones stated in \cite{BO}.

First, we formulate the properties of the process that 
directly follow from its construction
as a limit of finite Markov chains $T_n$
preserving measures $M_n$ on $\Kb_n$.
To prove them we use the convergence of measures $M_n$
(the Kingman theorem, \S\ref{s2.3})
and the properties of the up/down chains (\S\ref{s1.2}).

{\bf{}Invariant measure\/}
(Cf. \cite[Proposition 1.6]{BO}).
{\em{}The Poisson-Dirichlet distribution 
$\pd{\al}{\te}$ 
is an invariant measure for the process
$\x{\al}{\te}(t)$\/}.

This follows from the fact that each chain $T_n$ preserves the 
measure $M_n$.

{\bf{}Reversibility of the process\/}
(Cf. \cite[Proposition 1.7 and Theorem 7.3 (2)]{BO}).
{\em{}The process $\x{\al}{\te}(t)$ is reversible with respect to the measure $\pd{\al}{\te}$\/}.
	
This follows from the fact that each chain $T_n$ is reversible with respect
to $M_n$.

{\bf{}Convergence of finite-dimensional distributions\/}
(Cf. \cite[Propo\-sition 1.8]{BO}).
{\em{}Let $\x{\al}{\te}(t)$
and all the chains $T_n$ 
are viewed in equilibrium
(that is, starting from the invariant distribution). 
Then 
the finite-dimensional distributions for the $n$th chain converge, as $n\to\infty$,
to the corresponding finite-dimensional distributions of the process $\x{\al}{\te}(t)$. 
Here we
assume a natural scaling of time described before Proposition 
\ref{p4.3}\/}.

We proceed to the properties that
follow from the expression (\ref{f10}) for the 
pre-generator $A$ of the process $\x{\al}{\te}(t)$.

{\bf{} The spectrum of the Markov generator in $L^2\big(\sk,\pd{\al}{\te}\big)$.\/}
{\em{}The pre-generator $A$ acts in the space $\Fc\cong\Lao$.
Define an inner product in it:
\begin{equation*}
	(f,g)_{\pd{\al}{\te}}:=\int_{\sk}f(x)g(x)\,\pd{\al}{\te}(dx).
\end{equation*}
Then the space
$\Lao$ 
can be decomposed into the orthogonal direct sum
of eigenspaces of the operator $A$.
The spectrum of the operator $A$ looks as follows:
\begin{equation*}
	\left\{ 0 \right\}\cup\left\{ -\sigma_m\colon m=2,3,\dots \right\},\qquad 
	\sigma_m=m(m-1+\te).
\end{equation*}
The eigenvalue $0$ is simple
and the multiplicity of each $-\sigma_m$ 
is the number
of partitions of $m$
without parts equal to 1\/}.

The existence of such a decomposition of $\Lao$ follows from 
the fact that the pre-generator $A$ is symmetric (see the reversibility property above)
and preserves the filtration 
$\{\Fc^m\}$ of the space $\Fc\cong\Lao$
(see (\ref{f10})).
The fact that the operator $A$ is triangle
in the basis
$\{\m_\mu^\circ\}_{\mu\in\Kb,\;\kr{\mu}1=0}$
(compatible with the filtration $\{\Fc^m\}$)
implies the facts about the eigenstructure. 

{\bf{}The uniqueness of the invariant measure.\/}
{\em{}The measure $\pd{\al}{\te}$ 
is a unique invariant measure for the process
$\x{\al}{\te}(t)$\/}.

This can be proved similar to \cite[Theorem 7.3 (1)]{BO}.

{\bf{}Ergodicity.\/} 
{\em{}The process $\x{\al}{\te}(t)$ 
is ergodic 
with respect to the measure $\pd{\al}{\te}$, that is,
\begin{equation*}
	\lim_{t\to+\infty}\Big\|T(t)f-\int_{\sk}f(x)\,\pd{\al}{\te}(dx)\Big\|=0\quad\mbox{for all $f\in C(\sk)$},
\end{equation*}
where $\|\cdot\|$ is the supremum norm of the space $C(\sk)$\/}.

This follows from the existence of a spectral gap 
of the process' generator, see the eigenstructure above.
A detailed proof is given in \cite[Theorem 7.3 (3)]{BO}. 

\section{The pre-generator of the process in the simplex
as a differential operator}\label{s5}
This section performs a more detailed study of the
properties of the pre-generator $A\colon\Fc\to\Fc$
defined in Lemma \ref{p4.1}.
In \S\ref{s5.1} and \S\ref{s5.2} we assume the parameters $\al$ and $\te$ to
be of principal series, and
in \S\ref{s5.4} we study degenerate values of them.
\subsection{Pre-generator in natural coordinates}\label{s5.1}
Recall that the process $\x{\al}{\te}(t)$
in $\sk$ has the generator $\overline A$
(the closure of the pre-generator
$A\colon\Fc\to\Fc$).
The algebra $\Fc\cong\Lao$ is defined in \S\ref{s2},
and the operator $A$ is given by (\ref{f10})
in the basis
$\{\m_\mu^\circ\}_{\mu\in\Kb,\;\kr{\mu}1=0}$.
It turns out that in this basis the operator $A$ can be 
written as a second order differential operator 
in natural coordinates
$x_1,x_2,\dots$ on $\sk$.
First, let us express the functions
$\m_\mu^\circ(x)$ in these coordinates.
\begin{prop}\label{p5.1}
	For all
	$\mu\in\Kb$ such that $\kr{\mu}{1}=0$ the function
	$\m_\mu^\circ(x)$ on $\sk$ can be written as
	\begin{equation*}
		\m_\mu^\circ(x)=\sum_{i_1,\dots,i_{\ell(\mu)}}
		x_{i_1}^{\mu_1}\dots x_{i_{\ell(\mu)}}^{\mu_{\ell(\mu)}},\qquad x=(x_1,x_2,\dots)\in\sk,
	\end{equation*}
	where the sum is taken over all pairwise distinct indexes
	$i_1,\dots,i_{\ell(\mu)}$ from one to infinity.
\end{prop}
Note that if
$\kr{\mu}{1}>0$ and $\sum_{i=1}^{\infty}x_i<1$,
then this statement is false.
\begin{proof}
	The Proposition follows from a more general result by S.~Kerov
	\cite{Ke}. This result allows to evaluate
	any function of the form
	$\m_\la^\circ$, $\la\in\Kb$
	at any point $x\in\sk$.
	By $f(x_1,x_2,\dots)$
	denote the formal evaluation of the element $f\in\La$
	at the point $x\in\sk$.
	Let $\la\in\Kb$, $\kr{\la}{1}=r\ge0$. Then
	\begin{equation*}
		\m_\la^\circ(x)=\sum_{k=0}^{r}C_r^k\big(\gamma(x)\big)^k\m_{\la-k\cdot\square(1)}(x_1,x_2,\dots),\qquad x\in\sk.
	\end{equation*}
	Here $\m_\la^\circ(x)$ is the value of the function $\m_\la^\circ$ at the point $x$
	(this value was constructed in \S\ref{s2} using the factorization of the algebra $\La$);
	$\gamma(x)=1-\sum_{i=1}^{\infty}x_i$; and $\la-k\cdot\square(1)$ is
	the diagram obtained from $\la$ by deleting $k$ one-box rows.

	It is clear that the Proposition is a special case of this formula.
\end{proof}	
\begin{note}\label{p5.1.5}\rm{}
	It follows from the proof that 
	the function $g^{\circ}(x)$ on $\sk$
	for every element $g$ from the (non-factorized) algebra $\La$
	can be constructed 
	not only using the factorization of $\La$ (as explained in \S\ref{s2}).
	This can also be done directly.
	Namely, for every point 
	$x\in\nabla_\infty$ 
	(for which $\sum_{i=1}^{\infty}x_i=1$)
	we set $g^\circ(x)$ 
	to be the formal evaluation of the symmetric function $g\in\La$
	at the coordinates of the point
	$x=(x_1,x_2,\dots)$.
	We then extend the function
	$g^{\circ}(x)$ to the whole simplex $\sk$ by continuity.
	This method is suggested by the comment after the formula 
	\cite[(2.10)]{EK1}.
\end{note}

Now it is not hard to compute the pre-generator 
$A$ in natural coordinates.
Let $D$ denote the following formal expression:
\begin{equation}\label{f12}
	D=\sum_{i=1}^{\infty}x_i\frac{\partial^2}{\partial x_i^2}-
	\sum_{i,j=1}^{\infty}x_ix_j\frac{\partial^2}{\partial x_i\partial x_j}-
	\sum_{i=1}^{\infty}(\te x_i+\al)\frac{\partial}{\partial x_i}.
\end{equation}
It can be easily checked (using Proposition \ref{p5.1})
that for all
$\mu\in\Kb$ for which $\kr{\mu}1=0$ we have
\begin{equation*}
	D\m_\mu^\circ=A\m_\mu^\circ=-|\mu|(|\mu|-1+\te)\m_\mu^\circ+
	\sum_{\textstyle\genfrac{}{}{0pt}{}{c=1}{\mu_c\ge2}}^{\ell(\mu)}
	\mu_c(\mu_c-1-\al)\m^\circ_{\mu-\square(\mu_c)}.
\end{equation*}
It follows that
$Df=Af$ for all $f\in\Fc\subset C(\sk)$. Therefore,
the formula (\ref{f0.2}) from Introduction is justified.
\begin{note}\rm{}\label{p5.2}
	The RHS of (\ref{f12})
	can be understood in two equivalent ways as follows.
	Let $f\in\Fc$ and we want to compute
	$Df=Af$. Then we can either

	(1) express $f$ as a linear combination of vectors of the basis
	$\{\m_\mu^\circ\}_{\mu\in\Kb,\;\kr{\mu}1=0}$
	for the algebra $\Fc$,
	and then apply the operator $D=A$ to each $\m_\mu^\circ$ separately; or
	
	(2) compute $Af(x)$ first for $x\in\nabla_\infty$
	directly applying to $f$ the RHS of 
	(\ref{f12}), and then
	extend $Af(x)$ to the whole $\sk$ by continuity
	(cf. this method with Remark \ref{p5.1.5}).
\end{note}
\begin{note}\rm{}\label{p5.3}
	One might want to apply $D$
	(the RHS of (\ref{f12}))
	to functions $f$ not in $\Fc$.
	However, if for a function
	$f\in C(\sk)$ the expression
	$Df$ has a meaning, it is unevident
	that $f$
	enters the domain of the generator $\overline A$.
	It is even less evident that the expression 
	$Df$ coincides with $\overline Af$.

	Consider an example $f=x_1+\frac\al\te$, $\te\ne0$. It is clear that $Df=-\te f$.
	Had the equality
	$Df=\overline Af$ held, 
	then $f$ would be an eigenfunction of the generator 
	$\overline A$ corresponding to the eigenvalue
	$-\te$. Since
	$f\in L^2\big(\sk,\pd{\al}{\te}\big)$
	and $-\te$ 
	is not an eigenvalue of 
	$\overline A$ in this space
	(see \S\ref{s4.2}),
	we get a contraction.

	As shown in the paper \cite{Schm}, in the case $\al=0$
	the function $f=x_1$
	enters the domain of the Dirichlet form corresponding to the
	operator
	$\overline A$ 
	(and, therefore, the domain of the operator $\overline A^{1/2}$).
	Putting aside the question whether the function
	$f=x_1+\frac\al\te$
	(and, therefore, the coordinate function $x_1$) 
	enters the domain of $\overline A$,
	we conclude that the equality
	$Df=\overline Af$ is impossible. 
	
	This argument explains the following seeming paradox.
	If we treat the RHS of (\ref{f12})
	as a generator of a diffusion process
	(similarly to the finite-dimensional situation),
	then the drift vector
	at the point
	$x=(0,0,\dots)$ is equal to
	$(-\al,-\al,\dots)$ and therefore
	is directed ``away'' from the simplex $\sk$.
	But the action of the generator $\overline A$
	on the coordinate functions
	$x_1,x_2,\dots$ 
	is known not to be given by the RHS of
	(\ref{f12}). 
	It follows that in our situation the 
	coefficients of the first-degree derivatives
	in (\ref{f12}) cannot be interpreted as components of the drift vector.
\end{note}

\subsection{The pre-generator in moment coordinates}\label{s5.2}
First, recall the following definition.

Suppose we have a commutative algebra.
We say that every operator of multiplication
by an element of the algebra has zero order.
We say that an operator $Q$ has order $n$ if the commutator
of $Q$ with any operator of multiplication by an element of the algebra
has order $n-1$.
\begin{lemma}\label{p5.4}
	The pre-generator $A$ is a second order operator in the algebra $\Fc$.
\end{lemma}
\begin{proof}
	We must show that
	\begin{equation}\label{f19}
		\left[\big[[A,H_1],H_2\big],H_3\right]f(x)=0
	\end{equation}
	for all $h_1,h_2,h_3,f\in\Fc$ and $x\in\sk$, where by 
	$H_j$ we denote the operator of multiplication by $h_j$, $j=1,2,3$.

	From Remark \ref{p5.2} (2) it follows that 
	(\ref{f19}) holds for every $x\in\nabla_\infty$.
	Indeed, in this case the operator $A$ acts 
	according to the RHS of
	(\ref{f12}), and due to the linearity it
	remains to refer to the fact that
	$\left[ \left[ \left[ \partial^2/\partial x_i\partial x_j,H_1 \right],H_2 \right],H_3 \right]f(x)=0$
	and $\left[ \left[ \left[ \partial/\partial x_i,H_1 \right],H_2 \right],H_3 \right]f(x)=0$
	for all $1\le i,j<\infty$.

	For a general $x\in\sk$ the claim 
	(\ref{f19}) holds by continuity.
\end{proof}
Now,
the moment coordinates
$q_1,q_2,\dots$
are algebraically independent generators of the 
algebra $\Fc$ and the operator $A$ has second order in $\Fc$.
Therefore,
the action of $A$ on $\Fc$ is completely determined by 
the elements
$A1$, $Aq_i$, $A(q_iq_j)$, $i,j=1,2,\dots$.
It is clear that
$q_iq_j=q_{i+j+1}+\m^\circ_{(i+1,j+1)}$, $i\ge j$.
Using
(\ref{f10}), we find
\begin{equation*}
	\left.
	\begin{array}{lcl}
		\displaystyle
		A1&=& 0,\\\displaystyle Aq_i&=& -(i+1)(i+\te)q_i+(i+1)(i-\al)q_{i-1};\\\displaystyle
		A(q_iq_j)&=& q_iAq_j+q_jAq_i+2(i+1)(j+1)(q_{i+j}-q_iq_j),
	\end{array}
	\right.
\end{equation*}
where, by agreement, $q_0=1$.

Using the ``basis'' 
formal differential operators
$\partial/\partial q_i$ and $\partial^2/\partial q_i\partial q_j$
in the algebra $\Fc=\mathbb{R}\left[ q_1,q_2,\dots \right]$, we can write 
\begin{equation}\label{f15}\begin{array}{rcl}\displaystyle
	A&=&\displaystyle 
	\sum_{i,j=1}^{\infty}(i+1)(j+1)(q_{i+j}-q_iq_j)\frac{\partial^2}{\partial q_i\partial q_j}
	\\&&\displaystyle\qquad+
	\sum_{i=1}^{\infty}\left[ -(i+1)(i+\te)q_i+(i+1)(i-\al)q_{i-1} \right]\frac{\partial}{\partial q_i}.
	\end{array}
\end{equation}

This form of the operator $A$
allows to show that almost all sample paths of
the process $\x{\al}{\te}(t)$
are continuous.
This can be proved similarly to
\cite[Corollary 6.4 and Theorem 7.1]{BO}.

\subsection{Degenerate values of parameters}\label{s5.4}
Let us make mention what happens when the parameters $\al$ and $\te$
are of degenerate series (see \S\ref{s1.3} for the definition).
Consider the finite-dimensional ordered simplex
\begin{equation*}
	\nabla_N:=\left\{ (x_1,\dots,x_N)\in\mathbb{R^N}
	\colon x_1\ge\dots\ge x_N\ge0,\;\sum_{i=1}^{N}x_i=1 \right\},\quad N=2,3,\dots.
\end{equation*}
We can view $\nabla_N$ 
as a subset in 
$\nabla_\infty\subset\sk$.
The measure $\pd{\al}{\te}$ 
is concentrated on 
$\nabla_\infty\subset\sk$
for all values of $\al$ and $\te$ of principal or degenerate series.
\begin{prop}[{\cite[p. 62]{Pi3}}]\label{p5.5}
	When $\al=-\be<0$ and $\te=N\be$, the measure $\pd{-\be}{N\be}$ 
	is concentrated on $\nabla_N\subset\sk$.
        It coincides with the measure on $\nabla_N$
	that has the following
	density with respect to the Lebesgue measure:
	\begin{equation*}
		\pd{-\be}{N\be}(dx)=\frac{N!\,\Gamma(N\be)}{(\Gamma(\be))^{N}}\,
		x_1^{\be-1}\dots x_{N}^{\be-1}dx_1\dots dx_{N-1}.
	\end{equation*}

	Moreover, as $\be\to+0$ and $N\to+\infty$ such that $N\be\to\eta>0$,
	the following weak convergence of measures on $\sk$ holds:
	$\pd{-\be}{N\be}\to\pd{0}{\eta}$.
\end{prop}
It turns out that this result expands 
to the process $\x{\al}{\te}(t)$ in $\sk$.
First, we recall the diffusion processes considered
in the finite-dimensional simplex $\nabla_N$.

We are interested in processes that 
are called {\em{}approximating diffusions for the Wright-Fisher genetic 
model with symmetric mutation\/} in \cite{EK1}.
They were also studied in
\cite{Wa1,Gr} and other papers.
A detailed construction
of more general finite-dimensional diffusions
related to population genetics is explained in 
\cite[Chapter 10]{EK2}.

Let $\Fc_N=\mathbb{R}\left[ q_1,\dots,q_{N-1} \right]$ 
be a commutative unital algebra
freely generated by the moment coordinates on $\nabla_N$:
$q_k(x)=\sum_{i=1}^{N}x_i^{k+1}$, $k=1,\dots,N-1$.
It is clear that 
$\Fc_N\subset C(\nabla_N)$ is a dense subalgebra.
For each parameter $\eta>0$ 
we define an operator in $\Fc_N$:
\begin{equation}\label{f22}
	A_{N,\eta}:=\sum_{i=1}^{N}x_i\frac{\partial^2}{\partial x_i^2}-
	\sum_{i,j=1}^{N}x_ix_j\frac{\partial^2}{\partial x_i\partial x_j}+
	\frac{\eta}{N-1}\sum_{i=1}^{N}(1-Nx_i)\frac{\partial}{\partial x_i}.
\end{equation}
Note that this operator
(contrary to its infinite-dimensional analogue from \S\ref{s5.1})
can be defined by this formula on a wider subspace of 
$C(\nabla_N)$.
This subspace consists of twice continuously differentiable
under certain boundary conditions, see
\cite[(2.8)]{EK1}. 
Denote this extension by 
$\widetilde A_{N,\eta}$.
\begin{prop}[\cite{EK1}]\label{p5.6}
	{\rm{}(1)\/}
	The operator $\widetilde A_{N,\eta}$ (and, therefore, $A_{N,\eta}$)
	is closable in $C(\nabla_N)$. Denote by $\overline A_{N,\eta}$
	its closure (and, therefore, the closure of $A_{N,\eta}$);
	
	{\rm{}(2)\/} The closure $\overline A_{N,\eta}$ generates
        a diffusion process 	
	(that is, a strong Markov process with continuous sample paths
	which can start from any point and any probability distribution)
	in $\nabla_N$;

	{\rm{}(3)\/} 
	This process preserves the measure
	$\pd{-\frac{\eta}{N-1}}{\frac{N\eta}{N-1}}$
	defined above, and is reversible with respect to this measure.
\end{prop}
Denote this process by 
$\mathbf{Y}_{N,\eta}(t)$.
\begin{prop}[{\cite[Theorem 2.5]{EK1}}]\label{p5.7}
	As $N\to\infty$,
	the processes
	$\mathbf{Y}_{N,\eta}(t)$
	converge\footnote{In a certain sense explained
	in \cite[Theorem 2.5]{EK1}.}
	to the process $\x{0}{\eta}(t)$ constructed in \S\ref{s4}.
\end{prop}

It appears that the process
$\x{\al}{\te}(t)$
(which is constructed and studied throughout the present paper)
can be considered also for degenerate parameters.
First, we describe how the up/down Markov chains should be
modified in this case.

Let $\Kb(N)$ consist of all diagrams $\la\in\Kb$ such that
$\ell(\la)\le N$. 
It is again a graded graph.
Let $\al=-\be$, $\te=N\be$.
In this case each of the measures
$M_n$ defined in \S\ref{s1.3}
is positive everywhere on
$\Kb_n(N):=\Kb(N)\cap\Kb_n$.
It follows that we can consider the up/down Markov chains
on
$\Kb_n(N)$ similarly to \S\ref{s1.2}.
The place of the algebra $\La$ is taken by 
the algebra $\La_N$
of symmetric functions in $N$ variables.
It is clear that 
$\La_N/I\cong\Fc_N$, where $I=(p_1-1)\La$ is the ideal of the algebra
(see \S\ref{s2.1}).
Each embedding
$\iota_n$, $n=1,2,\dots$ (they were defined in \S\ref{s2.3})
now maps
$\Kb_n(N)$ into
$\nabla_N$.

\begin{prop}\label{p5.8}
	{\rm{}(1)\/}
	Let $T_n(N)$ 
	be the transition operator of the $n$th up/down Markov chain.
	In the same sense as in Lemma
	\ref{p4.1} there holds a convergence
	\begin{equation*}
		\lim_{n\to\infty}n^2(T_n(N)-{\bf1})f=A_{N,\be(N-1)}f\quad\mbox{for all $f\in\Fc_N$}.
	\end{equation*}

	{\rm{}(2)\/}
	The discrete semigroups
	$\left\{ {\bf1},T_n(N),T_n^2(N),\dots \right\}$ converge, as
	$n\to\infty$ (in the same sense as in Proposition \ref{p4.3}),
	to a continuous semigroup in the space
	$C(\nabla_N)$ generated by the operator
	$\overline A_{N,\be(N-1)}$.
\end{prop}
Thus, it is natural to say that
that the process
$\x{-\be}{N\be}(t)$ coincides with the process $\mathbf{Y}_{N,\be(N-1)}(t)$, that is,
$\x{-\be}{N\be}(t)$ is a finite-dimensional
diffusion process in
$\nabla_N$ with the pre-generator $A_{N,\be(N-1)}\colon\Fc_N\to\Fc_N$.
It follows that the diffusions in the finite-dimensional
simplexes studied in 
\cite{EK1} and many other works
arise as a special case of the two-parameter 
family of diffusions 
$\x{\al}{\te}(t)$ in $\sk$.

Proposition \ref{p5.7} 
can be interpreted as a convergence
(in the sense of \cite[Theorem 2.5]{EK1})
of processes
$\x{-\be}{N\be}(t)$ in $\nabla_N$ to the process
$\x{0}{\eta}(t)$ in $\sk$, as $\be\to+0$ and $N\to+\infty$ such that $N\be\to\eta>0$.
This convergence
strengthens
the second claim of Proposition
\ref{p5.5} concerning the convergence of invariant measures.

Kharkevich
Institute for Information Transmission Problems.
\end{document}